\documentclass[11pt]{amsart}
\usepackage{mathrsfs}
\usepackage{amssymb}

\usepackage{titletoc}
\usepackage{amscd}
\usepackage{amsmath, amssymb}
\usepackage{amsfonts}
\usepackage[colorlinks,linkcolor=blue,citecolor=blue, pdfstartview=FitH]{hyperref}

\numberwithin{equation}{section}

\theoremstyle{plain}
\newtheorem{thm}{Theorem}[section]
\newtheorem{theorem}[thm]{Theorem}
\newtheorem{lemma}[thm]{Lemma}
\newtheorem{corollary}[thm]{Corollary}
\newtheorem{proposition}[thm]{Proposition}

\theoremstyle{definition}

\newtheorem{question}[thm]{Question}

\newtheorem{remark}[thm]{Remark}

\newtheorem{definition}[thm]{Definition}

\newtheorem{example}[thm]{Example}

\newtheorem{defn-thm}[thm]{Definition-Theorem}
\newtheorem{conjecture}[thm]{Conjecture}

\newcommand{\sO}{{\mathcal O}}

\newcommand{\C}{{\mathbb C}}

\renewcommand{\P}{{\mathbb P}}

\newcommand{\T}{{\mathbb T}}

\renewcommand{\S}{{\mathbb S}}

\newcommand{\bp}{\bar{\partial}}

\newcommand{\ds}{\oplus}

\newcommand{\ts}{\otimes}

\newcommand{\btheorem}{\begin{theorem}}
\newcommand{\etheorem}{\end{theorem}}
\newcommand{\bproposition}{\begin{proposition}}
\newcommand{\eproposition}{\end{proposition}}
\newcommand{\bdefinition}{\begin{definition}}
\newcommand{\edefinition}{\end{definition}}
\newcommand{\bcorollary}{\begin{corollary}}
\newcommand{\ecorollary}{\end{corollary}}
\newcommand{\bproof}{\begin{proof}}
\newcommand{\eproof}{\end{proof}}
\newcommand{\bremark}{\begin{remark}}
\newcommand{\eremark}{\end{remark}}
\newcommand{\eexample}{\end{example}}
\newcommand{\bexample}{\begin{example}}
\newcommand{\la}{\langle}
\newcommand{\elemma}{\end{lemma}}
\newcommand{\blemma}{\begin{lemma}}
\newcommand{\ra}{\rangle}
\newcommand{\sq}{\sqrt{-1}}

\newcommand{\p}{\partial}

\renewcommand{\bar}{\overline}

\renewcommand{\phi}{\varphi}

\newcommand{\ee}{\end{eqnarray*}}
\newcommand{\be}{\begin{eqnarray*}}

\newcommand{\beq}{\begin{equation}}
\newcommand{\eeq}{\end{equation}}

\newcommand{\bd}{\begin{enumerate}}
\newcommand{\ed}{\end{enumerate}}

\renewcommand{\hat}{\widehat}
\renewcommand{\tilde}{\widetilde}

\renewcommand{\>}{\rightarrow}

\usepackage{fancyhdr}
\renewcommand{\leftmark}{{\scriptsize Hermitian manifolds with semi-positive holomorphic sectional curvature}}

\begin{document}
\title{Hermitian manifolds with semi-positive holomorphic sectional curvature}
\makeatletter
\let\uppercasenonmath\@gobble
\let\MakeUppercase\relax
\let\scshape\relax
\makeatother
\author{Xiaokui Yang}
\date{}

\address{ Current address of Xiaokui Yang: Morningside Center of Mathematics, and  Institute of Mathematics, and  Hua Loo-Keng Key Laboratory of
Mathematics, Academy of Mathematics \& Systems Science, The Chinese
Academy of Sciences, Beijing, China, 100190.}
\email{xkyang@amss.ac.cn}

\maketitle

\begin{abstract} We prove that a compact Hermitian manifold with
semi-positive but not identically zero holomorphic sectional
curvature has Kodaira dimension $-\infty$. As applications, we show
that Kodaira surfaces and hyperelliptic surfaces can not admit
Hermitian metrics with semi-positive holomorphic sectional curvature
although they have nef tangent bundles.

\end{abstract}

\section{Introduction}
In this note, we study compact Hermitian manifolds with
semi-positive \emph{holomorphic sectional curvature}. It is
well-known that, the holomorphic sectional curvature plays an
important role in differential geometry and algebraic geometry, e.g.
in  establishing the existence and nonexistence of rational curves
on projective manifolds. However, the relationships between
holomorphic sectional curvature and Ricci curvature, and the
algebraic positivity of the (anti-)canonical line bundles, and some
birational invariants  of the ambient manifolds are still
mysterious.
 In early 1990s, Yau proposed the following
question in his ``$100$ open problems in geometry"  (e.g.
\cite[Problem~67]{Yau} or \cite[p.392]{SY}):
\begin{question}\label{q} If $(X,\omega)$ is a compact K\"ahler manifold with
positive holomorphic sectional curvature, is $M$ unirational? Does
$X$ have negative Kodaira dimension?
\end{question}
\noindent At first, we answer  Yau's question  partially, but in a
more general setting.

 \btheorem\label{main} Let $(X,\omega)$ be a compact Hermitian
manifold with semi-positive holomorphic sectional curvature. If the
holomorphic sectional curvature is not identically zero, then $X$
has Kodaira dimension $-\infty$. In particular, if $(X,\omega)$ has
positive Hermitian holomorphic sectional curvature, then
$\kappa(X)=-\infty$. \etheorem

\noindent A complex manifold $X$ of complex dimension $n$ is called
complex parallelizable if there exist $n$  holomorphic vector fields
linearly independent everywhere. It is well-known that every complex
parallelizable manifold has flat curvature tensor and so identically
zero holomorphic sectional curvature. There are many non-K\"ahler
complex parallelizable manifolds of the form $G/H$, where $G$ is a
complex Lie group and $H$ is a discrete co-compact subgroup. On the
other hand, a compact complex parallelizable manifold is K\"ahler if
and only if it is a torus (\cite[Corollary~2]{Wang}). Hence, a
compact K\"ahler manifold $(X,\omega)$ with identically zero
holomorphic sectional curvature must be a torus.

As an application of Theorem \ref{main}, we obtain new examples of
K\"ahler and non-K\"ahler manifolds which can not support Hermitian
metrics with semi-positive holomorphic sectional curvature.

\bcorollary Let $X$ be a Kodaira surface or a hyperelliptic surface.
Then $X$ has nef tangent bundle, but $X$ does not admit a Hermitian
metric with semi-positive holomorphic {sectional} curvature.
\ecorollary

\noindent It is known  that  Kodaira surfaces and hyperelliptic
surfaces are complex manifolds with torsion anti-canonical line
bundles. Hence they are all complex Calabi-Yau manifolds. Note also
that Kodaira surfaces are all non-K\"ahler. Here, $X$ is said to be
a complex Calabi-Yau manifold if it has vanishing first Chern class,
i.e. $c_1(X)=0$. Moreover, \bcorollary\label{CY0} Let $X$ be a
compact Calabi-Yau manifold. If $X$ admits a K\"ahler metric with
semi-positive holomorphic sectional curvature, then $X$ is a torus.
\ecorollary

\noindent Note also that all diagonal Hopf manifolds are
non-K\"ahler Calabi-Yau manifolds with semi-positive holomorphic
sectional
curvature(\cite{LY14,Yang14}).\\

For reader's convenience, we present some recent progress about the
relationship between the positivity of holomorphic sectional
curvature and the algebraic positivity of the anti-canonical line
bundle, which is inspired by the following conjecture of Yau:

\begin{conjecture} \label{2222} If a compact K\"ahler manifold $(X,\omega)$   has strictly
negative holomorphic  sectional curvature, then the canonical line
bundle $K_X$ is ample.
\end{conjecture}
\noindent Bun Wong proved in \cite{Wong} that if $(X,\omega)$ is a
compact K\"ahler surface with negative holomorphic sectional
curvature, then the canonical line bundle $K_X$ is ample. Recently,
Heier-Lu-Wong showed in \cite{HLW10} that if $(X,\omega)$ is a
projective threefold with negative holomorphic sectional curvature,
then $K_X$ is ample. Moreover, by assuming the still open
``abundance conjecture" in algebraic geometry, they also confirmed
the conjecture for higher dimensional projective manifolds. On the
other hand, Wong-Wu-Yau proved in \cite{WWY} that if $(X,\omega)$ is
a compact projective manifold with Picard number $1$ and
quasi-negative holomorphic sectional curvature, then $K_X$ is ample.
As a breakthrough, Wu-Yau \cite{WY} confirmed  Conjecture \ref{2222}
when $X$ is projective. Building on their ideas, Tosatti and the
author proved  Conjecture \ref{2222} in full generality. More
precisely, we obtained

\btheorem[\cite{TY}]\label{TY} Let $(X,\omega)$ be a compact
K\"ahler manifold with nonpositive holomorphic sectional curvature.
Then the canonical line bundle $K_X$ is nef.  Moreover, if
$(X,\omega)$  has strictly negative holomorphic  sectional
curvature, then the canonical line bundle $K_X$ is ample. \etheorem

 \noindent For more related discussions on this topic, we refer to
\cite{Wong,HLW10,WWY,HLW14,WY,TY,HW15} and the references therein.
One may also wonder whether similar statements hold for compact
K\"ahler manifolds with positive holomorphic sectional curvature.
However,

\bexample  Let $Y$ be the Hirzebruch surface
$Y=\P\left(\sO_{\P^1}(-k)\ds \sO_{\P^1}\right)$ for $k\geq 2$. It is
is proved (\cite{Hitchin75} or \cite[p.292]{SY}) that $Y$ has a
smooth K\"ahler metric with positive holomorphic sectional
curvature. But the anti-canonical line bundle $K_Y^{-1}$ is not
ample although $K_Y^{-1}$ is known to be effective.   For more
details, see Example \ref{ex}.\eexample

\noindent As an important structure theorem, Heier and Wong proved
in
  \cite{HW12} that projective manifolds with positive total scalar curvature are uniruled. In partiuclar, projective manifolds with
  positive holomorphic sectional curvature are uniruled.\\

\section{Preliminaries}

Let $(E,h)$ be a Hermitian holomorphic vector bundle over a compact
complex manifold $X$ with Chern connection $\nabla$. Let
$\{z^i\}_{i=1}^n$ be the  local holomorphic coordinates
  on $X$ and  $\{e_\alpha\}_{\alpha=1}^r$ be a local frame
 of $E$. The curvature tensor $R^\nabla\in \Gamma(X,\Lambda^2T^*X\ts E^*\ts E)$ has components \beq R_{i\bar j\alpha\bar\beta}= -\frac{\p^2
h_{\alpha\bar \beta}}{\p z^i\p\bar z^j}+h^{\gamma\bar
\delta}\frac{\p h_{\alpha \bar \delta}}{\p z^i}\frac{\p
h_{\gamma\bar\beta}}{\p \bar z^j}\eeq (Here and henceforth we
sometimes adopt the Einstein convention for summation.)
 In
particular, if $(X,\omega_g)$ is a  Hermitian manifold,
$(T^{1,0}M,\omega_g)$ has  Chern curvature components \beq R_{i\bar
j k\bar \ell}=-\frac{\p^2g_{k\bar \ell}}{\p z^i\p\bar z^j}+g^{p\bar
q}\frac{\p g_{k\bar q}}{\p z^i}\frac{\p g_{p\bar \ell}}{\p\bar
z^j}.\eeq  The (first) Chern-Ricci form $Ric(\omega_g)$ of
$(X,\omega_g)$ has components
$$R_{i\bar j}=g^{k\bar \ell}R_{i\bar jk\bar \ell}=-\frac{\p^2\log\det(g)}{\p z^i\p\bar z^j}$$
and it is well-known that the Chern-Ricci form represents the first
Chern class of the complex manifold $X$ (up to a factor $2\pi$). The
Chern scalar curvature $s$ of $(X,\omega_g)$ is defined as \beq
s=g^{i\bar j} R_{i\bar j}. \eeq For a Hermitian manifold $(X,
\omega_g)$, we define the torsion tensor \beq T_{ij}^k=g^{k\bar
\ell}\left(\frac{\p g_{j\bar \ell}}{\p z^i}-\frac{\p g_{i\bar
\ell}}{\p z^j}\right).\eeq By using elementary Bochner formulas
(e.g. \cite[Lemma~3.3]{LY14}, or \cite[Lemma~A.6]{LY12}), we have
\beq \bp^*\omega=-\sq T_{ki}^kdz^i.\label{bp*}\eeq Indeed, we have
 $[\bp^*,L]=\sq \left(\p+\tau\right)$
with $\tau=[\Lambda,\p\omega]$. When it acts on constant $1$, we
obtain \be \bp^*\omega=\sq\tau(1)&=&\sq
\Lambda\left(\p\omega\right)\\
&=&-\sq g^{k\bar\ell}\left(\frac{\p g_{i\bar\ell}}{\p z^k}-\frac{\p
g_{k\bar\ell}}{\p z^i}\right)dz^i\\
&=&-\sq T_{ki}^kdz^i.\ee

 Let $(X,\omega)$ be a compact Hermitian manifold.
$(X,\omega)$ has positive (resp. semi-positive)  holomorphic
sectional curvature, if for any nonzero vector $\xi=(\xi^1,\cdots,
\xi^n)$,
$$R_{i\bar j k\bar \ell}\xi^i\bar\xi^j\xi^k\bar\xi^\ell>0\ \ \ \text{(resp. $\geq 0$)}.$$

\section{Hermitian manifolds with semi-positive holomorphic sectional
curvature}\label{section3}

In this section, we discuss the relationship between the holomorphic
sectional curvature and the Kodaira dimension of the ambient
manifold. It is also well-known that, on Hermitian
 manifolds, there are many curvature notations and the curvature relations are more
 complicated than the relations in the K\"ahler case because of the non-vanishing of
 the torsion tensor
 (e.g.\cite{LY12,LY14}).

\btheorem\label{main1}  Let $(X,\omega)$ be a compact Hermitian
manifold with semi-positive holomorphic sectional curvature. If the
holomorphic sectional curvature is not identically zero, then $X$
has Kodaira dimension $\kappa(X)=-\infty$. \etheorem

\bproof  At a given point $p\in X$, the maximum holomorphic
sectional curvature is defined to be
$$ H_p:=\max_{W\in T^{1,0}_pX,|W|=1}H(W),$$ where $H(W):=R(W,\bar W,
W,\bar W)$. Since $X$ is of finite dimension, the maximum can be
attained. Suppose the holomorphic sectional curvature is not
identically zero, i.e.  $ H_p>0$ for some $p\in X$. For any $q\in
X$. We assume $g_{i\bar j}(q)=\delta_{ij}$. If $\dim_{\C} X=n$ and
$[\xi^1,\cdots,\xi^n]$ are the homogeneous coordinates on
$\P^{n-1}$, and $\omega_{FS}$ is the Fubini-Study metric of
$\P^{n-1}$. At point $q$, we have the following well-known
identity(e.g. \cite[Lemma~4.1]{LYJAG}): \beq\int_{\P^{n-1}}R_{i\bar
j k\bar \ell}\frac{\xi^i\bar
\xi^j\xi^k\bar\xi^\ell}{|\xi|^4}\omega^{n-1}_{FS}=R_{i\bar j k\bar
\ell}\cdot
\frac{\delta_{ij}\delta_{k\ell}+\delta_{i\ell}\delta_{kj}}{n(n+1)}=\frac{s+\hat
s}{n(n+1)}.\label{scalar}\eeq where $ s$ is the Chern scalar
curvature of $\omega$ and $\hat s$ is defined as \beq \hat
s=g^{i\bar \ell} g^{k\bar j}R_{i\bar j k\bar \ell}.\eeq Hence if
$(X,\omega)$ has semi-positive holomorphic sectional curvature, then
$s+\hat s$ is a non-negative function on $X$. On the other hand, at
point $p\in X$,  $s+\hat s$ is strictly positive. Indeed, since
$H_p>0$, there exists a nonzero vector $\xi\in T_p^{1,0}X$ such that
$H(\xi)=R_{i\bar j k\bar \ell}\frac{\xi^i\bar
\xi^j\xi^k\bar\xi^\ell}{|\xi|^4}>0$. By (\ref{scalar}), the
integrand is quasi-positive over $\P^{n-1}$, and so $s+\hat s$ is
strictly positive at $p\in X$. Note that in general if $(X,\omega)$
is not K\"ahler, $s$ and $\hat s$ are not the same.  By
\cite[Section ~4]{LY14},  we have the relation \beq s=\hat s+\la
\bp\bp^*\omega,\omega\ra.\label{2}\eeq Indeed, we compute \be s-\hat
s&=&g^{i\bar j}g^{k\bar \ell}\left(R_{i\bar j k\bar \ell}-R_{k\bar j
i\bar \ell}\right)=g^{i\bar j}g^{k\bar \ell}\left(\nabla_{\bar
j}\left(\frac{\p g_{i\bar \ell}}{\p z^k}-\frac{\p g_{k\bar \ell}}{\p
z^i}\right)\right)\\
&=&g^{i\bar j}\nabla_{\bar j}T_{ki}^k=g^{i\bar j}\frac{\p
T_{ki}^k}{\p\bar z^j}=\la\bp\bp^*\omega,\omega\ra\ee where we use
 formula (\ref{bp*}) in the last identity. Therefore, we have
\beq \int_X \hat s \omega^n=\int_X s\omega^n-\int_X
|\bp^*\omega|^2\omega^n.\label{3}\eeq

Next we use Gauduchon's conformal method (\cite{Ga2,Ga3}, see also
\cite{Balas1,Balas2}) to find a Hermitian metric $\tilde \omega$ in
the conformal class of $\omega$ such that $\tilde \omega$ has
positive
 Chern scalar curvature $\tilde s$.

 Let
$\omega_G=f_0^{\frac{1}{n-1}}\omega$ be a Gauduchon metric ( i.e.
$\p\bp\omega_G^{n-1}=0$ ) in the conformal class of $\omega$ for
some strictly positive weight function $f_0\in C^\infty(X)$
(\cite{Ga2,Ga3}). Let $s_G,\hat s_G$ be the corresponding scalar
curvatures with respect to the Gauduchon metric $\omega_G$. Then we
have
\begin{eqnarray} \int_Xs_G\omega_G^n\nonumber&=&-n\int_X \sq
\p\bp\log\det(\omega_G)\wedge \omega_G^{n-1}\\
\nonumber&=&-n\int_X\left(\sq\p\bp\log\det(\omega)+\frac{n}{n-1}\sq\p\bp\log
f_0\right)\wedge \omega_G^{n-1}\\
\nonumber&=&-n\int_X\sq\p\bp\log\det(\omega)\wedge \omega_G^{n-1}\\
\nonumber&=&-n\int_X f_0\sq\p\bp\log\det(\omega)\wedge
\omega^{n-1}\\
&=&\int_X f_0 s\omega^n,\label{4} \end{eqnarray} where we use the
Stokes' theorem and the fact that $\omega_G$ is Gauduchon in the
third identity. Similarly, by using the proof of formula
(\ref{bp*}), we have the relation \be
\bp^*_G\omega_G&=&\sq\Lambda_G(\p\omega_G)\\
&=&\sq\Lambda(\p\omega)-\sq f_0^{-\frac{1}{n-1}}\left(\frac{\p
f_0^{\frac{1}{n-1}}}{\p z^k}-n\frac{\p f_0^{\frac{1}{n-1}}}{\p
z^k}\right)dz^k
\\&=&\bp^*\omega+\sq\p
\log f_0.\ee Since $\omega_G$ is Gauduchon,  we obtain
 \begin{eqnarray} \int_X \la \bp\bp_G^*\omega_G,\omega_G\ra\omega_G^n
\nonumber&=&n\int_X\bp\bp_G^*\omega_G\wedge \omega_G^{n-1}=n\int_X
\bp\bp^*\omega\wedge \omega_G^{n-1}\\&=&\int_X
f_0\la\bp\bp^*\omega,\omega\ra\omega^n.\label{5} \end{eqnarray}
 By using a similar equation as (\ref{3}) for $s_G,\hat s_G$ and $\omega_G$, we obtain \begin{eqnarray} \int_X
\hat s_G\omega_G^n\nonumber&=&\int_X
s_G\omega_G^n-\int_X\la\bp\bp^*_G\omega_G,\omega_G\ra\omega_G^n\\&=&\int_X
f_0 s\omega^n-\int_Xf_0\la\bp\bp^*\omega,\omega\ra\omega^n=\int_X
f_0 \hat s\omega^n.\label{6}\end{eqnarray} where we use equations
(\ref{4}), (\ref{5}) in the second identity, and (\ref{2}) in the
third identity.
 Therefore, if $s+\hat s$ is quasi-positive,
we obtain \begin{eqnarray} \int_X
s_G\omega_G^n\nonumber&=&\frac{\int_X (s_G+\hat
s_G)\omega_G^n}{2}+\frac{\int_X (s_G-\hat
s_G)\omega_G^n}{2}\\&=&\frac{\int_X (s_G+\hat
s_G)\omega_G^n}{2}+\frac{\|\bp^*_G\omega_G\|^2}{2}=\frac{\int_X
f_0(s+\hat s)\omega^n}{2}+\frac{\|\bp^*_G\omega_G\|^2}{2}>0
\label{balanced}\end{eqnarray} where the third equation follows from
(\ref{4}) and (\ref{6}).

Next, there exists a Hermitian metric $h$ on $K^{-1}_X$ which is
conformal to $\det(\omega_G)$ on $K^{-1}_X$ such that the scalar
curvature $s_h$ of $(K_X^{-1},h)$ {with respect to $\omega_G$ is a
constant}, and more precisely we have
$$\displaystyle{s_h=-
\text{tr}_{\omega_G}\sq\p\bp\log h=\frac{\int_X
s_G\omega_G^n}{\int_X\omega^n_G}}.$$ Indeed, let $f\in C^\infty(X)$
be a strictly positive function satisfying \beq s_G-
\text{tr}_{\omega_G}\sq\p\bp f=\frac{\int_X
s_G\omega_G^n}{\int_X\omega^n_G} \label{hopf}\eeq then
$h=f\det(\omega_G)$ is the metric we need. Note that the existence
of  solutions to (\ref{hopf}) is well-known by Hopf's lemma.

Finally, we deduce that the conformal metric $$ \tilde \omega:=
f^{\frac{1}{n}}f_0^{\frac{1}{n-1}}\omega=f^{\frac{1}{n}}\omega_G$$
is a Hermitian metric with positive Chern scalar curvature. Indeed,
the Chern scalar curvature $\tilde s$  is, \be \tilde
s&=&-\text{tr}_{\tilde\omega }\sq \p\bp\log\det(\tilde
\omega^n)\\&=&-\text{tr}_{\tilde \omega}\sq\p\bp\log h\\
&=&-f^{-\frac{1}{n}} \text{tr}_{\omega_G}\sq\p\bp\log h \\
&=&f^{-\frac{1}{n}}\frac{\int_X
s_G\omega_G^n}{\int_X\omega^n_G}>0.\ee

Hence, if $\sigma\in H^0(X,mK_X)$ for some positive integer $m$, by
the standard Bochner formula with respect to the metric $\tilde
\omega$, one has
 \beq
tr_{\tilde\omega}\sq\p\bp|\sigma|^2_{\tilde\omega}=|\nabla'
\sigma|_{\tilde \omega}^2+m\tilde s\cdot|\sigma|_{\tilde\omega}^2
\eeq where $|\bullet|_{\tilde\omega}$ is the pointwise norm on
$mK_X$ induced by $\tilde\omega$ and $\nabla'$ is the $(1,0)$
component of the Chern connection on $mK_X^{}$. Since $\tilde s$ is
strictly positive, by maximum principle we have
$|\sigma|^2_{\tilde\omega}=0$, i.e. $\sigma=0$. Now we deduce the
Kodaira dimension of $X$ is $-\infty$.
 \eproof

It is  easy to see that, on a K\"ahler manifold
 $(X,\omega)$, if the total scalar curvature $\int_X s\omega^n$ is
 positive, then $\kappa(X)=-\infty$ (e.g. \cite[Theorem~1]{Yau2}
or \cite[Theorem~1.1]{HW12}). However, in general, it is not true
for
 non-K\"ahler metrics
  which can be seen from the following
 example.

 \bexample Let $(\T^2,\omega)$ be a torus with the flat metric.
 For any non-constant real smooth function $f\in C^\infty(\T^2)$, the Hermitian metric
 $\omega_f=e^f\omega$ has strictly positive total Chern scalar
 curvature and $\kappa(\T^2)=0$. Indeed,
 $\det(\omega_f)=e^{2f}\det(\omega)$ and
 $$Ric(\omega_f)=-\sq\p\bp\log\det(\omega_f)=Ric(\omega)-2\sq\p\bp f=-2\sq\p\bp f.$$
 The total scalar curvature of $\omega_f$ is given by
 \be \int s_f\cdot \omega_f^2&=&\int tr_{\omega_f} Ric(\omega_f) \cdot\omega^2_f\\&=&2\int
 Ric(\omega_f)\wedge \omega_f=-4\int \sq \p\bp f\wedge e^f \omega\\
 &=&4\int \left(\sq \p f\wedge \bp f\right) e^f \omega\\
 &=&4\|\p f\|^2_{\omega_f}>0
 \ee
 since $f$ is not a constant function, where we use the Stokes' theorem in the
 fourth identity.
 \eexample

 Note that, a special case of Theorem \ref{main1} is proved
in \cite{Balas2} that when $X$ is a surface or a threefold and
$(X,\omega)$ has strictly positive holomorphic sectional curvature,
then $X$ has Kodaira dimension $-\infty$.

As an application of Theorem \ref{main1}, we have
\bcorollary\label{example} Let $X$ be a Kodaira surface or a
hyperelliptic surface. Then $X$ has nef tangent bundle, but $X$ does
not admit a Hermitian metric with semi-positive holomorphic
\emph{sectional} curvature. \ecorollary \bproof It is well-known
that the holomorphic tangent bundles of Kodaira surfaces or
hyperelliptic surfaces are nef (e.g. \cite{DPS} or \cite{Yang14}).
On the other hand, if $X$ is either a Kodaira surface or a
hyperelliptic surface, then $X$ has torsion canonical line bundle,
i.e. $K_X^{\ts m}=\sO_X$ for some positive integer $m$
(\cite[p.244]{BHPV}). In particular, we have $\kappa(X)=0$. Suppose
$X$ has a Hermitian metric $\omega$ with semi-positive holomorphic
sectional curvature, by Theorem \ref{main1}, $(X,\omega)$ has
constant zero holomorphic sectional curvature. Then $(X,\omega)$ is
a K\"ahler surface \cite[Theorem~1]{BG}. Since all Kodaira surfaces
are non-K\"ahler, we deduce that Kodaira surfaces can not admit
Hermitian metrics with semi-positive holomorphic sectional
curvature. Suppose $(X,\omega)$ is a hyperelliptic surface with
constant zero holomorphic sectional curvature. So $\omega$ is a
K\"ahler metric with constant zero holomorphic sectional curvature,
and we deduce $(X,\omega)$ is flat since the curvature tensor is
determined by the holomorphic sectional curvature. Indeed, for any
$Y,Z\in T^{1,0}_pX$, expand
$$R(Y+\lambda Z, \bar {Y+\lambda Z}, Y+\lambda Z, \bar{Y+\lambda Z})\equiv0$$
into powers of $\lambda$ and $\bar \lambda.$ Using the K\"ahler
symmetry, the $|\lambda|^2$ term gives $R(Y,\bar Y, Z,\bar Z)=0$.
Now if we expand
$$R(Y+\lambda Z, \bar{Y+\lambda Z}, A+\mu B, \bar{A+\mu B})\equiv0$$
into powers of $\lambda, \bar\lambda, \mu, \bar \mu$, the
$\bar{\lambda \mu}$ term gives $R(Y,\bar Z, A,\bar B)=0$ for any
$(1,0)$-vectors  $Y,Z,A,B\in T^{1,0}_pX$. Since $(X,\omega)$ is
flat,  $X$ is a complex parallelizable manifold (e.g.
\cite[Proposition~2.4]{DLV} and \cite{AS}). However, it is proved in
\cite[Corollary~2]{Wang} that a complex parallelizable manifold is
K\"ahler if and only if it is a torus. This is a contradiction.
\eproof

Let $X$ be a  complex manifold. $X$ is said to be a complex
Calabi-Yau manifold if $c_1(X)=0$. \bcorollary\label{CY} Let $X$ be
a compact Calabi-Yau manifold. If $X$ admits a K\"ahler metric with
semi-positive holomorphic sectional curvature, then $X$ is a torus.
\ecorollary \bproof Let $X$ be a compact K\"ahler Calabi-Yau
manifold, then it is well-known  that (e.g. \cite[Theorem~1.5]{T}),
$K_X$ is a holomorphic torsion, i.e. there exists a positive integer
$m$ such that $K_X^{\ts m}=\sO_X$. In particular, $\kappa(X)=0$.
Suppose $X$ has a smooth K\"ahler metric $\omega$ with semi-positive
holomorphic sectional curvature, then by Theorem \ref{main1}, $X$
has constant zero holomorphic sectional curvature. As shown in
Corollary \ref{example}, $(X,\omega)$ is flat and so it is a complex
parallelizable manifold, i.e. $X$ is a torus. \eproof

\bremark As shown in \cite{Yang14}, the Hopf surface $H_{a,b}$ (and
every diagonal Hopf manifold \cite{LY14}) has a Hermitian metric
with semi-positive holomorphic bisectional curvature. Since
$b_2(H_{a,b})=b_2(\S^1\times \S^3)=0$, we see $c_1(H_{a,b})=0$ and
so $H_{a,b}$ is a non-K\"ahler Calabi-Yau manifold with
semi-positive holomorphic sectional curvature. \eremark

Finally, we want to use the following well-known example to
demonstrate that the positivity of the holomorphic sectional
curvature can not imply the ampleness of the anti-canonical line
bundle although the negativity of the holomorphic sectional
curvature does imply the ampleness of the canonical line bundle
(e.g. \cite{TY} or Theorem \ref{TY}).

\bexample\label{ex}  Let $Y$ be the Hirzebruch surface
$Y:=\P\left(\sO_{\P^1}(-k)\ds \sO_{\P^1}\right)$ for $k\geq 2$ which
is a $\P^1$-bundle over $\P^1$. It is known (\cite{Hitchin75} or
\cite[p.292]{SY}, or \cite{AAH}) that $Y$ has a smooth K\"ahler
metric with positive holomorphic sectional curvature. Next, we show
$K_Y^{-1}$ is not ample although $K_Y^{-1}$ is known to be
effective.
  Let $E:=\sO_{\P^1}(k)\ds
\sO_{\P^1}$, $Y=\P(E^*)$ and $\sO_Y(1)$ be the tautological line
bundle of $Y$. The following adjunction formula is well known
(e.g.\cite[p.89]{Laza1}) \beq K_Y=\sO_{Y}(-2)\ts \pi^*(K_{\P^1}\ts
\det E)\eeq where $\pi$ is the projection $Y=\P(E^*)\>\P^1$. In
particular, we have \beq \sO_Y(2)= K^{-1}_Y\ts
\pi^*(\sO_{\P^1}(k-2)).\eeq Suppose $K_Y^{-1}$ is ample, then
$\sO_Y(1)$ is aslo ample since $k\geq 2$. Therefore, by definition
(\cite{Ha}), $E=\sO_{\P^1}(k)\ds \sO_{\P^1}$ is an ample vector
bundle. This is a contradiction. \eexample

 \section*{Acknowledgements} The author would like to thank
 K.-F. Liu,  V. Tosatti, B. Weinkove and S.-T. Yau for many
valuable discussions.  The author would also like to thank the
anonymous referees whose comments and suggestions helped improve and
clarify the paper.

\end{document}